\documentclass[12pt,a4paper]{article}
\usepackage{amsmath,amsfonts,amssymb,amscd}

\usepackage[a4paper]{geometry}

\usepackage{ulem}
\usepackage{xcolor}
\usepackage{graphics}

\usepackage{hyperref}

\def\Ker{\operatorname{Ker}}

\def\ad{\operatorname{ad}}

\def\d{\operatorname{d}}

\def\Im{\operatorname{Im}}

\def\Lie{\operatorname{Lie}}

\newcounter{th}
\def\t{\refstepcounter{th}{\bf \noindent{Theorem} \arabic{th}. }}

\newcounter{prop}
\def\prop{\refstepcounter{prop}{\bf \noindent{Proposition} \arabic{prop}. }}

\newcounter{lem}
\def\lem{\refstepcounter{lem}{\bf \noindent{Lemma} \arabic{lem}. }}

\newcounter{de}
\def\de{\refstepcounter{de}{\bf \noindent{Definition} \arabic{de}. }}

\newcounter{ex}
\def\ex{\refstepcounter{ex}{\bf \noindent{Example} \arabic{ex}. }}

\begin{document}

\begin{center}
    {\LARGE{\bf Commutative $n$-ary superalgebras with an invariant skew-symmetric form\footnote{ Supported by  AFR-grant, University of Luxembourg.
     }}}
\end{center}

\begin{center}
     E.G. Vishnyakova
\end{center}

\bigskip

\begin{abstract}
We study $n$-ary commutative superalgebras and $L_{\infty}$-algebras that possess a skew-symmetric invariant form, using the derived bracket formalism.
 This class of superalgebras includes for instance Lie algebras and their $n$-ary generalizations, commutative associative and Jordan algebras with an invariant form.  
 We give a classification of anti-commutative $m$-dimensional $(m-3)$-ary algebras with an invariant form, and a classification of real simple $m$-dimensional Lie $(m-3)$-algebras with a positive definite invariant form up to isometry. Furthermore,
 we develop the Hodge Theory for $L_{\infty}$-algebras with a symmetric invariant form, and we describe quasi-Frobenius structures on skew-symmetric $n$-ary algebras.

\end{abstract}

\bigskip

\section{Introduction}

\noindent \texttt{Derived bracket formalism.} The derived bracket
approach was successfully used in different areas of mathematics: in
Poisson geometry, in the theory of Lie algebroids and Courant
algebroids, BRST formalism, in the theory of Loday algebras and different
types of Drinfeld Doubles. For detailed introduction we
recommend a beautiful survey of Y.~Kosmann-Schwarzbach
\cite{Kosman1}. 

The idea of the formalism is the following. 
{\it One fixes an algebra $L$, usually a Lie superalgebra, and
constructs another multiplication on the same vector space (or some
subspace) using derivations of $L$ and the (iterated) multiplication in
$L$. One obtains a class of new algebras, which properties can be
studied using original algebra $L$.} 
For example, using this formalism we can obtain all Poisson structures on a manifold $M$ from the canonical Poisson algebra on $T^*M$
 as was shown by Th.~Voronov in \cite{Ted Poisson}. Voronov's idea
allows A.~Cattaneo and M.~Zambon \cite{Catt} to introduce a unified approach to the reduction of Poisson manifolds.
Another example was suggested in \cite{Ted Der Bracket homotopy} and \cite{Ted Der Bracket all algebras}, where a series of strongly homotopy algebras was obtained from a given Lie superalgebra.

We use this formalism to study $n$-ary commutative superalgebras with an invariant skew-symmetric form. More precisely,
 consider a vector superspace $V$ with a non-degenerate even skew-symmetric form $(\,,)$. There exists a natural Lie superalgebra structure on $S^*(V)$, where $S^*(V)$ is the symmetric power of $V$.
  The main observation is that {\it we get all commutative $n$-ary and
  strongly homotopy superalgebras on $V$ with the invariant skew-symmetric
   form $(\,,)$. In other words, the property of these $n$-ary superalgebras having an invariant skew-symmetric form is encoded by the Lie superalgebra $S^*(V)$.} 
The observation that using the superalgebra $S^*(V)$ we can obtain all Lie algebras with an invariant
symmetric form was made by B.~Kostant and S.~Sternberg in \cite{KosSter}. The superalgebra $S^*(V)$ was also
used in Poisson Geometry to study for instance Lie bialgebras and Drinfeld Doubles, see \cite{Kosman1,Kosman2}, \cite{LecRoger} and others.

\medskip

\noindent\texttt{Multiple generalizations of Lie algebras.} Using the
derived bracket formalism we can study $n$-ary commutative superalgebras with a skew-symmetric invariant form. This class
of superalgebras includes for instance different $n$-ary generalizations of a Lie algebra with a symmetric invariant form.
First of all let us give a short review of such generalizations.

Multiple generalizations arise usually from different readings of the Jacobi identity. For example, the Jacobi identity
for a Lie algebra is equivalent to the statement that all adjoint operators  are derivations of this Lie algebra. If we
use this point  of view for the $n$-ary case we come to the notion of a {\it Filippov $n$-algebra} \cite{Filippov}.
V.T.~Filippov considered alternating n-ary algebras $A$ satisfying
the following Jacobi identity:
\begin{equation}\label{eq_ Jacobi Filippov}
\{ a_1, \ldots ,a_{n-1}, \{ b_1, \ldots, b_n\}\} =
\sum
\{ b_1, \ldots, b_{i-1}\{ a_1, \ldots, a_{n-1}, b_i\}, \ldots, b_n\},
\end{equation}
where $a_i,\,b_j\in A$. In other words, the operators $\{ a_1, \ldots ,a_{n-1}, - \}$ are
derivations of the $n$-ary bracket $\{ b_1, \ldots, b_n \}$. Such algebras appear naturally in  Nambu
mechanics \cite{Nambu} in the contecst of Nambu-Poisson manifolds, in supersymmetric
gravity theory and in supersymmetric gauge theories, the Bagger-Lambert-Gustavsson Theory, see \cite{n-ary algebras: a review}.

Another natural n-ary generalization of the Jacobi identity has the following form:
\begin{equation}\label{eq_ Jacobi SH}
\sum (-1)^{(I,J)} \{ \{ a_{i_1},\ldots, a_{i_n}\}, a_{j_1}, \ldots, a_{j_{n-1}}\} = 0,
\end{equation}
where the sum is taken over all ordered unshuffle multi-indexes
$I=(i_1, \ldots,i_n)$ and $J= (j_1, \ldots, j_{n-1})$  such that
$(I,J)$ is a permutation of $(1,\ldots, 2n-1)$. We will call such algebras {\it Lie $n$-algebras}. 
 This type of $n$-ary algebras was considered for instance by P.~Michor and A.~Vinogradov in  \cite{MV} and by P.~Hanlon
and M.L.~Wachs \cite{HW}.
The homotopy case was studied in \cite{SS} in context of the
Schlesinger-Stasheff homotopy algebras and  $L_{\infty}$-algebras. Such algebras are related to
the Batalin-Fradkin-Vilkovisky theory and to the string field theory, see \cite{Lads Stasheff}.
In \cite{VV1} A.M.~Vinogradov and M.M.~Vinogradov
proposed a three-parameter family of n-ary algebras  such that for some $n$
the above discussed  structures appear as particular cases.

The theory of Filippov $n$-algebras is relatively well-developed. For instance,
there is a classification of simple real and complex Filippov $n$-algebras and an analog of the Levi decomposition \cite{Ling}. W.X.~Ling in \cite{Ling} proved that there exists only one simple finite-dimensional Filippov  $n$-algebra over an algebraically closed field of characteristic $0$ for any $n>2$.
The simple Filippov $n$-superalgebras in the finite and infinite dimensional
case were studied in \cite{Kac}. It was shown there that there are no simple linearly
compact Filippov $n$-superalgebras which are not Filippov $n$-algebras, if $n > 2$, and a
 classification of linearly compact Filippov $n$-algebras was given.

In this paper we give a classification of anti-commutative $m$-dimensional $(m-3)$-ary algebras with a symmetric invariant form  over $\mathbb R$ and $\mathbb C$ up to isometry in terms of coadjoint orbits of the Lie group $\mathrm{SO}(V)$. In the real positive definite case we give  a classification of simple algebras of this type. Our result can be formulated as follows: {\it almost all real anti-commutative $m$-dimensional $(m-3)$-ary algebras with a symmetric invariant positive definite form are simple. The exceptional cases are: the trivial $(m-3)$-ary algebra and the $(m-3)$-ary algebras that corresponds to decomposible element.}
We also give a classification of real (simple) $m$-dimensional Lie $(m-3)$-algebras with a symmetric invariant positive definite form.

\medskip

\noindent\texttt{Hodge decomposition for real strongly homotopy algebras.} A definition of a
strongly homotopy Lie algebra (or $L_{\infty}$-algebra or sh-algebra) was given by Lada and Stasheff in \cite{Lads Stasheff}.
For more about strongly homotopy algebras see also  \cite{Lada Markl}, \cite{Ted Der Bracket homotopy}, \cite{Ted Der Bracket all algebras}.
Another result of our paper is a Hodge Decomposition for  real metric pure odd strongly homotopy algebras.
An observation here is that we can obtain easily such kind of decomposition using the derived bracket formalism.

We can also use this formalism to define the Hodge operator on a Riemannian compact oriented manifold $M$.
Indeed, in this case there exists the metric on cotangent space $T^*M$ that is induced by Riemannian metric on
the tangent space $TM$. Then we can define a Poisson bracket on $\bigwedge T^*M$, see \cite{Roytenberg}, and
repeat the construction of the Hodge operator given in the present paper.

\medskip

\noindent \texttt{Quasi-Frobenius structures.} We conclude our paper with a description of quasi-Frobenius
structures on anti-commutative $n$-ary algebras. Our result is as follows. {\it Assume that $n$ is even. There is  a one-to-one
correspondence between quasi-Frobenius structures on an anti-commutative $n$-ary algebra and maximal isotropic subalgebras in $T^*_0$-extension on this  algebra.}

\section{Commutative $n$-ary superalgebras with an invariant skew-sym\-metric form }

\subsection{Main definitions}

Let $V = V_{\bar {0}} \oplus V_{\bar {1}}$ be a finite dimensional $\mathbb{Z}_2$-graded vector space
 over the field $\mathbb{K}$, where $\mathbb{K} = \mathbb{R}$ or $\mathbb{C}$.  If $a\in V$ is a
 homogeneous element, we denote by $\bar{a}\in \mathbb{Z}_2$ the parity of $a$. As
 usual we assume that elements in $\mathbb{K}$ are even. Recall that a bilinear form $(\,,)$ on $V$ is called
 {\it even} (or {\it odd}) if the corresponding  linear map $V\otimes V \to \mathbb{K}$
 is even (or odd).  A bilinear form is called {\it skew-symmetric} if
$
(a,b)= -(-1)^{\bar{a}\bar{b}}(b,a)
$
 for any homogeneous elements $a,b\in V$.

\medskip

\de \label{def invarian mult, symmetric}
 $\bullet$ An {\it $n$-ary superalgebra structure} on $V$ is an $n$-linear map 
 \begin{align*}
V\times \cdots \times & V \longrightarrow V,\\
(a_1,\ldots,a_n)& \mapsto \{a_1,\ldots,a_n\}.
 \end{align*}

$\bullet$ An $n$-ary superalgebra structure is called {\it commutative} if
\begin{equation}\label{eq symmetry}
\{a_1,\ldots, a_i, a_{i+1}, \ldots, a_n \} = (-1)^{\bar{a_i}\bar{a}_{i+1}} \{a_1,\ldots, a_{i+1}, a_i, \ldots, a_n\}
\end{equation}
for any homogeneous $a_i, a_{i+1}\in V$.

\smallskip

$\bullet$ A commutative $n$-ary  superalgebra structure is called {\it invariant with respect to the form} $(\,,)$ if  the following holds:
\begin{equation} \label{eq quadratic}
(a_0, \{a_{1}, \ldots, a_{n}\}) = (-1)^{\bar{a}_0 \bar{a}_1}(a_1, \{a_0,a_2, \ldots, a_{n}\})
\end{equation}
for any homogeneous $a_{i}\in V$.

\medskip

We will write {\it a commutative invariant $n$-ary superalgebra structure} or {\it a commutative invariant $n$-ary superalgebra} as a shorthand for {\it a commutative $n$-ary superalgebra structure on $V$ that is invariant with respect to the form $(\,,)$}.

\medskip

\ex \label{ex symmertic graded invariant n -ary superalgebras} The class of commutative invariant $n$-ary superalgebras includes for instance the following algebras.

$\bullet$ {\it Anti-commutative algebras on $V = V_{\bar 1}$ with an invariant symmetric form.}  Indeed, in this case the conditions
(\ref{eq symmetry}) and (\ref{eq quadratic}) are equivalent to the following conditions:
\begin{equation} \label{ex skew-symmetric invariant}
\{a,b\} = -\{b,a\}, \quad (\{a,b\}, c) = (a,\{b,c\}).
\end{equation}
In particular, all Lie algebras with an invariant symmetric form are of this type.

\smallskip

$\bullet$ {\it Commutative algebras on $V = V_{\bar 0}$ with an invariant skew-symmetric form.}  In this case from (\ref{eq symmetry})
and (\ref{eq quadratic}) it follows:
\begin{equation} \label{ex symmetric invariant}
\{a,b\} = \{b,a\}, \quad (\{a,b\}, c) = -(a,\{b,c\}).
\end{equation}
In particular, commutative associative and Jordan algebras with an invariant skew-symmetric form are of this type.

\smallskip

$\bullet$ {\it Anti-commutative $n$-ary algebras on $V = V_{\bar 1}$ with an invariant symmetric form.} In this case 
the condition (\ref{eq quadratic}) is equivalent to the following condition:
$$
(y,\{x_1,\ldots,x_{n-1},z\}) = (-1)^n (\{y,x_1,\ldots,x_{n-1}\},z)
$$
that is more familiar for physicists. In particular, anti-commutative $n$-ary algebras satisfying (\ref{eq_ Jacobi Filippov}) with an invariant symmetric form are of this type. Such algebras are used in  the Bagger-Lambert-Gustavsson model (BLG-model), see \cite{n-ary algebras: a review} for details.

\medskip

\noindent {\bf Remark.} For a commutative algebra usually one
considers  the following invariance condition: $(\{a,b\}, c) =
(a,\{b,c\})$. If in addition we assume that  the form $(\,,)$ is
skew-symmetric and non-degenerate, we obtain $2(ab,c)=0$ for all $a,b,c \in V$, therefore $ab=0$. In our case we do not have such
additional restrictive relations.

\subsection{Derived bracket and commutative invariant $n$-ary superalgebras}

Let  $V$ be as above. We
denote by $S^n V$ the $n$-th symmetric power of $V$ and we put $S^*
V = \bigoplus\limits_n S^n V$. The superspace $S^* V$ possesses a natural
structure $[\,,]$ of a Poisson superalgebra. It is defined by the
following formulas:
$$
[x,y]: = (x,y), \quad x,y\in V;
$$
$$
[v,w_1\cdot w_2]: = [v,w_1]\cdot w_2 + (-1)^{vw_1} w_1\cdot [v,w_2],
$$
$$
[v,w] = -(-1)^{vw}[w,v],
$$
where $v,w, w_i$ are homogeneous elements in $S^* V$. One can show that the multiplication $[\,,]$ satisfies the  graded Jacobi identity:
$$
[v,[w_1, w_2]] = [[v, w_1], w_2] + (-1)^{\bar{v} \bar{w_1}}[w_1,[v,w_2] ].
$$
This Poisson superalgebra is well-defined. Indeed, we can repeat the
argument from \cite[Page 65]{KosSter}  for vector superspaces. The
idea is to show that this superalgebra is induced by the
Clifford superalgebra corresponding to $V$ and $(\,,)$.

Let us take any element $\mu \in S^{n+1} V$. Then we can define an n-ary superalgebra structure on $V$ in the following way:
\begin{equation}\label{eq der product}
\{a_{1}, \ldots, a_{n}\}: = [a_1,[\ldots,[a_n,\mu]\ldots]], \,\,\, a_i\in V.
\end{equation}
We will denote the corresponding superalgebra by $(V,\mu)$ and we will call the element $\mu$ the {\it derived
potential} of $(V,\mu)$.  
The $n$-ary superalgebras of type $(V,\mu)$ have the following two properties:

\begin{itemize}
\item The multiplication (\ref{eq der product}) is {\it commutative}. (This was noticed in \cite{Ted Der Bracket homotopy}.)
Indeed, using Jacobi identity for  $S^* V$ we have:
 \begin{align*}
[a_1,[a_2,\ldots,[a_n,\mu]\ldots]] =& [[a_1,a_2],[\ldots,[a_n,\mu]\ldots]] +\\
(-1)^{\bar{a}_1\bar{a}_{2}} [a_2,[a_1,\ldots&,[a_n,\mu]\ldots]] = (-1)^{\bar{a}_1\bar{a}_{2}} [a_2,[a_1,\ldots,[a_n,\mu]\ldots]].
\end{align*}
We used the fact that $[[a_1,a_2],[\ldots,[a_n,\mu]\ldots]]=0$, because $[a_1,a_2]\in \mathbb{K}$.
 Similarly we can prove  the commutativity relation for other $a_i$.

\item The $n$-ary superalgebra structure (\ref{eq der product}) is invariant. 
 Indeed,
  \begin{align*}
 (a_0, \{a_{1}, \ldots, a_{n}\}) & =  [a_0, [a_1,[a_2,\ldots,[a_n,\mu]\ldots]]] = \\
 (-1)^{\bar{a}_0\bar{a}_1}&[a_1, [a_0,[a_2,\ldots,[a_n,\mu]\ldots]]] = 
  (-1)^{\bar{a}_0\bar{a}_1}(a_1, \{a_0,a_2, \ldots, a_{n}\}).
\end{align*}

\end{itemize}

We conclude this section with the following  observation.

\medskip

\prop\label{Main obseravation graded} {\it Assume that $V$ is finite dimensional and $(\,,)$ is
non-degenerate.  Any commutative invariant $n$-ary superalgebra structures can be obtained by construction
(\ref{eq der product}). }

\medskip

\noindent{\it Proof.} Denote by  $\mathcal{A}_n$  the vector space
of commutative invariant $n$-ary superalgebra  structures on $V$ and by $\mathcal{L}_{n+1}$ the vector space of symmetric $(n+1)$-linear maps from $V$ to $\mathbb{K}$. Clearly,  $\dim \mathcal{L}_{n+1} = \dim S^{n+1}V$. Since $(\,,)$ is non-degenerate, Formula (\ref{eq der
product}) defines an injective linear map $S^{n+1}V \to
\mathcal{A}_n$. 
 We can also define an injective linear map $\mathcal{A}_n \to \mathcal{L}_{n+1}$ in the following way:
$$
\mathcal{A}_n\ni \mu \longmapsto L_{\mu} \in \mathcal{L}_{n+1}, \quad L_{\mu}(a_1,\ldots,a_{n+1}) = (a_1,\mu(a_2,\ldots,a_{n+1})).
$$
Note that $L_{\mu}$ is symmetric since $\mu$ defines an invariant superalgebra structure.   Summing up, we  have the following sequence
of injective maps or isomorphisms:
$$
S^{n+1}V \hookrightarrow \mathcal{A}_n  \hookrightarrow \mathcal{L}_{n+1} \simeq S^{n+1}V.
$$
Since $V$ is finite dimensional, we get $S^{n+1}V \simeq \mathcal{A}_n$.$\Box$

\section{Examples of commutative invariant $n$-ary superalgebras}

Usually one studies  superalgebras with an invariant form in the
following way.  One considers for example a Lie algebra or a Jordan
algebra  and assumes that the multiplication
in the algebra satisfies the following additional condition: it is
invariant with respect to a non-degenerate (skew)-symmetric form.
The derived bracket formalism permits to express for instance Jacobi, Filippov and Jordan
identities in terms of derived potentials and the Poisson bracket on $S^* V$. In this case
the additional invariance condition is fulfilled automatically.

\subsection{Strongly homotopy Lie algebras with an invariant\\ skew-symmetric form}

We follow Th.~Voronov \cite{Ted Der Bracket homotopy} in conventions
concerning $L_{\infty}$-algebras. 
 We set $I^k:= (i_1, \ldots,
 i_k)$ and $J^l:= (j_1, \ldots,
  j_l)$, where  $i_1<\cdots <i_k$ and $j_1<\cdots <j_l$. We denote 
 $$
  a_{I^k}:= (a_{i_1},\ldots, a_{i_k}),\quad a_{J^l}:= (a_{j_1},\ldots, a_{j_l})\quad  \text{and} \quad a^s:= (a_1,\ldots, a_s),
  $$
   where $a_s\in V$. We put $[a_{I^k},\mu]:= [a_{i_1}, \ldots
[a_{i_k},\mu]]$ and $[a^{s},\mu]:= [a_{1}, \ldots
[a_{s},\mu]]$, where $\mu \in S^{*}V$.

\medskip

\de \label{de_L infinity algebra} A vector superspace $V$ with a
sequence of odd commutative $n$-linear  maps $\mu_n $, where $n \geq 0$, is
called an {\it $L_{\infty}$-algebra} if
the following generalized Jacobi identities hold:
\begin{equation}\label{generalized Jacobi}
\sum\limits_{k+l=n}\,\, \sum\limits_{(I^k,J^l)} (-1)^{(I^k,J^l)}  \mu_{l+1}(a_{I^l}, \mu_{k}( a_{J^k} ) ) =0, \quad n \geq 0.
\end{equation}
Here $(I^k,J^l)$ is a unshuffle permutation of $(1,\ldots, n)$ and $(-1)^{(I^k,J^l)}$ is the sign obtained using the sign rule for the permutation $(I^k,J^l)$ of homogeneous elements $a_1,\ldots, a_n\in V$.

\medskip

\de \label{de_L infinity algebra quadratic} An
$L_{\infty}$-algebra structure $(\mu_n)_{n\geq 0}$ on $V$ is called {\it invariant}  if all
$\mu_n$ are invariant in  the sense of
Definition \ref{def invarian mult, symmetric}.

\medskip

The following statement can be deduced from \cite[Theorem $1$]{Ted Der Bracket homotopy} and Proposition \ref{Main obseravation graded}. For completeness we give here a proof in our notations and agreements. 

\medskip

\prop\label{prop quadratic L_infinity algebra} {\it Invariant $L_{\infty}$-algebra structures  on $V$ are in one-to-one correspondence with odd
elements $\mu\in S^*(V)$ such that $[\mu,\mu]\in S^0V= \mathbb K$.}

\medskip

\noindent{\it Proof.} 
Our objective is to show that  $[\mu,\mu]\in \mathbb K$ is equivalent to (\ref{generalized Jacobi}) together with the invariance condition. 
Let us take any odd element $\mu=\sum\limits_k \mu_k\in
S^*V$, where $\mu_k\in S^{k+1}V$.  The condition $[\mu,\mu]\in \mathbb K$ is
equivalent to the following conditions
$$
\sum\limits_{k+l=n}  [\mu_l,\mu_k]=0, \quad  n \geq 1.
$$
Note that $[\mu_0,\mu_0]$ is always an element is $\mathbb K$.
This is equivalent to:
$$
\sum\limits_{k+l=n} [a^{n-1}, [\mu_l,\mu_k]]=0, \quad n \geq 1.
$$
 Furthermore, we have:
\begin{align*}
[a^{n-1},[\mu_l,\mu_k]] =  &\sum\limits_{(I^l,J^{k-1})} (-1)^{(I^l,J^{k-1})+ \bar{a}_{J^{k-1}}}  [[a_{I^l},\mu_l], [a_{J^{k-1}},\mu_{k}]] + 
\\
&\sum\limits_{(I^{l-1},J^{k})} (-1)^{(I^{l-1},J^{k})+ \bar{a}_{J^k}} [[a_{I^{l-1}},\mu_l], [a_{J^k},\mu_{k}]],
\end{align*}
where $\bar{a}_{J^{k-1}}$ and $\bar{a}_{J^k}$ are the parities of $a_{J^{k-1}}$ and $a_{J^k}$, respectively.
Denote by $\tilde{\mu}_s$ the $s$-linear map defined by $\tilde{\mu}_s(a_1,\ldots,a_s):= [a^s,\mu_s]$. We get:
\begin{align*}
[[a_{I^l},\mu_l], [a_{J^{k-1}},\mu_{k}]]&= \tilde\mu_{k}(\tilde\mu_l(a_{I^l}), a_{J^{k-1}});\\
[[a_{I^{l-1}},\mu_l], [a_{J^k},\mu_{k}]]& = (-1)^{(\bar{a}_{J^{l-1}}-1)(\bar{a}_{J^{k}}-1)+1} \tilde\mu_{l}(\tilde\mu_k(a_{I^k}), a_{J^{l-1}}).
\end{align*}
Further,
\begin{align*}
[a^{n-1},[\mu_l,\mu_k]] =
\sum\limits_{(I^{l},J^{k-1})} (-1)^{(I^{l},J^{k-1}) + \bar{a}_{J^{k-1}}} \tilde\mu_{k}(\tilde\mu_l(a_{I^l}), a_{J^{k-1}})  +\\
\sum\limits_{(I^{k},J^{l-1})} (-1)^{(J^{k},I^{l-1})+  \bar{a}_{J^{l-1}}} \tilde\mu_{l}(\tilde\mu_k(a_{I^k}), a_{J^{l-1}})=\\
\sum\limits_{(J^{k-1},I^l)} (-1)^{(J^{k-1},I^l)} \tilde\mu_{k}(a_{J^{k-1}},\tilde\mu_l(a_{I^l})) +\\
\sum\limits_{(I^{l-1},J^k)} (-1)^{(I^{l-1},J^k)} \tilde\mu_{l}(a_{I^{l-1}},\tilde\mu_k(a_{J^k})).
\end{align*}
Using the equolities:
\begin{align*}
\sum_{k+l=n}\sum\limits_{(J^{k-1},I^l)} (-1)^{(J^{k-1},I^l)} &\tilde\mu_{k}(a_{J^{k-1}},\tilde\mu_l(a_{I^l})) =\\
& \sum_{k'+l=n-1}\sum\limits_{(J^{k'},I^l)} (-1)^{(J^{k'},I^l)} \tilde\mu_{k'+1}(a_{J^{k'}},\tilde\mu_l(a_{I^l}));\\
\sum_{k+l=n}\sum\limits_{(I^{l-1},J^k)} (-1)^{(I^{l-1},J^k)} &\tilde\mu_{l}(a_{I^{l-1}},\tilde\mu_k(a_{J^k}))=\\
&\sum_{k+l'=n-1}\sum\limits_{(I^{l'},J^k)} (-1)^{(I^{l'},J^k)} \tilde\mu_{l'+1}(a_{I^{l'}},\tilde\mu_k(a_{J^k}))
\end{align*}
we see that
\begin{equation}\label{eq_[a^n-1,mu] = n-1 Jacobi}
[a^{n-1},\sum\limits_{k+l=n} [\mu_l,\mu_k]]=2 \sum_{k'+l=n-1}\sum\limits_{(J^{k'},I^l)} (-1)^{(J^{k'},I^l)} \tilde\mu_{k'+1}(a_{J^{k'}},\tilde\mu_l(a_{I^l})).
\end{equation}
Therefore, $[a^{n-1},\sum\limits_{k+l=n} [\mu_l,\mu_k]]=0$ is equivalent to the generalized $(n-1)$-Jacobi identity for the invariant $L_{\infty}$-algebra $\{\tilde{\mu}_s\}$. Conversely,  if an invariant $L_{\infty}$-algebra $\{\tilde{\mu}_s\}$ is given, its derived potential $\mu= \sum\limits_s\mu_s$, where $\mu_s\in S^{s+1}V$ corresponds to $\tilde{\mu}_s$ by Proposition \ref{Main obseravation graded}, must satisfy the condition $[\mu,\mu]\in \mathbb K$.
$\Box$

\medskip

\noindent {\bf Corollary.} {\it Assume that $V=V_{\bar 1}$ and $n$ is even. Anti-commutative invariant $n$-ary algebra structures, where $n>0$, on $V$ satisfying the Jacobi identity (\ref{eq_ Jacobi SH})  are in one-to-one correspondence with
elements $\mu\in S^{n+1}(V)$ such that $[\mu,\mu]=0$.}

\medskip

\noindent {\it Proof.} In this case the equation \ref{eq_[a^n-1,mu] = n-1 Jacobi} has the form:
$$
[a^{2n-1},[\mu,\mu]] = 2\sum\limits_{(I,J)} (-1)^{(I,J)} \tilde\mu(a^{I},\tilde\mu(a^{J})).
$$
Here $I=(i_1,\ldots, i_{n-1})$ and $J=(j_1,\ldots, j_{n})$ such that $i_1<\cdots< i_{n-1}$, $j_1<\cdots< j_{n}$ and $I\cup J = \{1,\ldots, 2n-1\}$. Since $n$ is even we have:
$$
\sum\limits_{(I,J)} (-1)^{(I,J)} \tilde\mu(a^{I},\tilde\mu(a^{J})) = - \sum\limits_{(J,I)} (-1)^{(J,\,I)} \tilde\mu(\tilde\mu(a^{J}), a^{I}),
$$
and the equality
$$
\sum\limits_{(J,I)} (-1)^{(J,\,I)} \tilde\mu(\tilde\mu(a^{J}), a^{I}) =0
$$
is equivalent to the Jacobi identity (\ref{eq_ Jacobi SH}) for $n$-ary algebra structure $\tilde\mu$.$\Box$

\medskip

\subsection{Filippov algebras with an invariant symmetric form}

The class of commutative symmetric superalgebras includes $n$-ary algebras introduced by Filippov \cite{Filippov}. 

\medskip

\de \label{de_Filippov algebra} $\bullet$ An anti-commutative  $n$-ary algebra is
called a {\it Filippov $n$-algebra} if it satisfies the Jacobi identity (\ref{eq_ Jacobi Filippov}).

$\bullet$ We say that a Filippov $n$-algebra $V$
 has an {\it invariant symmetric form} $(\,,)$ if the following holds 
$$
(y,\{x_1,\ldots, x_{n-1},z\}) = (-1)^n (\{y,x_1,\ldots, x_{n-1}\},z)
$$
for any $x_i,y,z\in V$. 

\medskip

Filippov $n$-algebras with an invariant form are described in the following proposition. The idea of the proof we borrow from \cite{VV1}.

\medskip

\prop\label{prop quadratic Filippov algebra} {\it
Assume that $V= V_{\bar 1}$. Invariant Filippov $n$-algebra structures on $V$ are in one-to-one correspondence with elements $\mu \in  S^{n+1}V$ such that the following equation holds:
$$
[\mu_{a^{n-1}}, \mu] =0
$$
for all $a^{n-1} = (a_1,\ldots, a_{n-1})$ and $\mu_{a^{n-1}}:= [a_1,\ldots, [a_{n-1},\mu]]$.

}

\medskip

\noindent{\it Proof.} We need to show that $[\mu_{a^{n-1}}, \mu] =0$ is equivalent to (\ref{eq_ Jacobi Filippov}) together with the invariance condition. Let us take $b_1\ldots, b_{n}\in V$. We have:
$$
\begin{array}{c}
[\mu_{a^{n-1}}, [b_1,\ldots [b_n,\mu]] ] = \sum\limits_{i=1}^n [b_1,\ldots [ [\mu_{a^{n-1}},b_i],\ldots  [b_n,\mu] ]] + 
[b_1,\ldots [b_n,[\mu_{a^{n-1}}, \mu]]] ]
\end{array}
$$
Further, using (\ref{eq der product}), we get:
 \begin{align*}
[\mu_{a^{n-1}}, [b_1,\ldots [b_n,\mu]] ] =  -\{\{b_1,\ldots, b_{n} \}, a_1\ldots,& a_{n-1}\} =\\
&(-1)^n\{ a_1,\ldots, a_{n-1}, \{b_1,\ldots, b_{n} \}\};
\end{align*}
and
 \begin{align*}
[b_1,\ldots [ [\mu_{a^{n-1}},b_i],\ldots  [b_n,\mu] ]] = & -\{b_1\ldots,b_{i-1}, \{ b_{i},a_1\ldots, a_{n-1}\}, b_{i+1},\ldots, b_{n}\} =\\
&(-1)^n \{b_1\ldots,b_{i-1}, \{a_1\ldots, a_{n-1},  b_{i}\}, b_{i+1},\ldots, b_{n}\};
\end{align*}
 Hence, we have:
$$
\begin{array}{rr}
\{ a_1,\ldots, a_{n-1}, \{b_1,\ldots, b_{n} \}\} = &\sum\limits_{i=1}^n \{b_1\ldots,b_{i-1}, \{a_1\ldots, a_{n-1},  b_{i}\}, b_{i+1},\ldots, b_{n}\} + \\
&(-1)^n  [b_1,\ldots [b_n,[\mu_{a^{n-1}}, \mu]]] .
\end{array}
$$
Therefore, the condition $[b_1,\ldots [b_n,[\mu_{a^{n-1}}, \mu]]] = 0$ is equivalent to 
(\ref{eq_ Jacobi Filippov}) together with the invariance condition. 
The proof is complete.$\Box$

\subsection{Jordan algebras with a  skew-symmetric invariant form}

First of all let us recall the definition of a Jordan algebra.

\medskip

\de  A {\it Jordan algebra} is a commutative algebra  such that the multiplication satisfies the following condition:
$$
(xy)(xx) = x(y(xx)).
$$

\medskip

\de  We say that a Jordan algebra $V$ has an {\it invariant skew-symmetric form} $(\,,)$ if the following holds: 
$$
(ab,c) = -(a,bc)
$$
for any $a,b,c\in V$. 
\medskip

A description of invariant Jordan algebrs structures on $V$ is given in the following proposition. 

\medskip

\prop\label{prop symplectic Jordan algebras}
 {\it Let $V$ be a pure even vector space with a non-degenerate skew-symmetric form $(\,,)$. Invariant Jordan algebra structures on $V$ are in one-to-one correspondence with elements 
  $A \in S^{3}V$   such that the following identity holds:
 $$
 [A_x, A_{[A_x,x]}] =0,
 $$
 where $A_x = [x,A]$.
}

\medskip

\noindent {\it Proof.} By Proposition \ref{Main obseravation graded}
any invariant commutative algebra structure $(x,y) \mapsto xy$ on $V$  can be obtained by the derived bracket construction. Denote by
$A$ its derived potential. In other words, we have:
$$
x y= [x,[y,A]].
$$
Further,
$$\begin{array}{c}
(xy)(xx) =  [[y,A_x], [[x,A_x], A ]];\quad

x(y(xx)) = - [A_x, [y, [[x,A_x], A]]].
\end{array}
$$
Using the Jacobi identity for the Poisson algebra $S^*V$, we get:
$$
\begin{array}{c}
[A_x, [y, [[x,A_x], A]]] = [[A_x,y], [[x, A_x], A]] + [y, [A_x, [[x,A_x], A] ] ].
\end{array}
$$
We see that this equation is equivalent to
$$
- x(y(xx)) = - (xy)(xx) +  [y, [A_x, [[x,A_x], A] ] ].
$$
Hence, the algebra $(V,A)$ is Jordan if and only if
$$
 [y, [A_x, [[x,A_x], A] ] ]=0
$$
for all $x,y\in V$. The last condition is equivalent to
$$
[A_x, [[x,A_x], A] ] = 0
$$
for all $x\in V$.$\Box$

\subsection{Associative algebras with a skew-symmetric invariant form}

\medskip

\prop\label{prop symplectic associative algebras}
{\it Let $V$ be a pure even vector space with a non-degenerate skew-symmetric form $(\,,)$. Invariant associative algebra structures on $V$ are in one-to-one correspondence with elements 
 $\mu \in  S^{3}V$ such that the following identity holds:
$$
[\mu_a,\mu_b] = 0
$$
 for all $a,b\in V$. Here  $\mu_x = [x,\mu]$ for $x\in V$.
}

\medskip

\noindent {\it Proof.}
  Let us use the notation:
 $
 a\circ b:= [a,[b,\mu]].
 $
   We have to show that the associativity relation for $\circ$ together with the invariance condition is equivalent to $[\mu_a,\mu_b] = 0$ for all $a,b\in V$.
Indeed,
 \begin{align*}
a\circ (b \circ c)= [a,[[b,[c,\mu]]\mu]] = - [a,[\mu,[b,[c,\mu]]] = 
- [\mu_a, [b,&[c,\mu]]] =\\ -[[\mu_a,b], [c,\mu]] - [ b, [\mu_a,\mu_c]] = [[b,\mu_a]&, [c,\mu]]  - [ b, [\mu_a,\mu_c]]=\\
 &(b\circ a)\circ c  - [ b, [\mu_a,\mu_c]].
 \end{align*}
Therefore, the associativity relation for $\circ$ together with the invariance condition and  the equality $[\mu_a,\mu_c] = 0$ for all $a,c\in V$ are equivalent.$\Box$

\section{Hodge operator and its applications}

\subsection{$*$-operator and $n$-ary algebras}

 Let $V$ be a pure odd vector space of dimension  $m$ with a non-degenerate
  skew-symmetric even bilinear form $(\,,)$. Recall that means that $(a,b) = (b,a)$ for all $a,b\in V$.
  Let us choose a normalized orthogonal basis $(e_i)$ of $V$. Denote by $L:= e_1  \dots e_m$ the top form
   corresponding to the chosen basis.  We define the operator $*:S^{p} V \to S^{m-p} V$ by the following formula:
 \begin{equation}\label{eq * operator coordinates}
  *(x_{1}  \dots x_{p}) = [x_{1}, [ \dots[x_{p} , L]]].
 \end{equation}
  In particular, we have:
$$
*(e_{i_1}  \dots e_{i_p}) = [e_{i_1}, [ \dots[e_{i_p} , L]]] = (-1)^{\sigma} e_{j_1}  \dots e_{j_{m-p}},
$$
where $\sigma(1, \ldots, m) = (i_p,\ldots, i_1,j_1,\ldots, j_{m-p})$.
 Clearly, this definition depends only on orientation of $V$ and on the bilinear form $(\,,)$. Note that $*: S^pV \to S^{m-p}V$
  is an isomorphism for all $p$. This follows for example from the following formula:
$$
**(e_{i_1}  \dots e_{i_p}) = (-1)^{\frac{m(m-1)}{2}} e_{i_1}  \dots e_{i_p}.
$$
The following well-known result we can easily prove using derived bracket formalism:

\medskip

\prop\label{prop lin operators preserving the form} {\it The vector
space $\mathfrak{so}(V)$ of  linear operators preserving the form
$(\,,)$ is isomorphic to $S^2(V)$.}

\medskip

\noindent{\it Proof.} The isomorphism is given by the formula $w \longmapsto \ad w,$
where $w\in S^2(V)$ and $\ad w (v):= [w,v]$ for $v\in V$. Indeed, for all $v_1, v_2\in V$ we have:
\begin{align*}
0 = \ad w([v_1,v_2])  = [w, [v_1,v_2]] = [[w,v_1],v_2] + [v_1,&[w,v_2]] =\\
& ([w,v_1],v_2)  + (v_1,[w,v_2]).
 \end{align*}
 Clearly, this map is injective. We complete the proof observing that the dimensions of $\mathfrak{so}(V)$ and $S^2(V)$ are equal.$\Box$

\medskip

We have seen in previous sections that elements from $S^{n+1}V$
corresponds to invariant $n$-ary algebra structures on $V$. The
existence of the $*$-operator for $V=V_{\bar 1}$ leads to the idea
that $n$-ary and $(m-n)$-ary algebras can have some common
properties. In particular such algebras have the same algebra of
orthogonal derivations.

\medskip

\de A {\it derivation} of an $n$-ary algebra $(V,\mu)$ is a linear map $D:V\to V$ such that
$$
D(\{v_1,\ldots, v_n\}) = \sum_j \{v_1,\ldots, D(v_j),\ldots, v_n\}.
$$
We denote by $\mathrm{IDer} (\mu)$ the vector space of all derivations of the algebra $(V,\mu)$ preserving the form $(\,,)$.

\medskip

\prop\label{prop derivations} {\it Let us take any $w\in S^2(V)$ and $\mu \in S^{n+1}(V)$.

\noindent{\bf a.} We have:
$$
\mathrm{IDer} (\mu) =  \{w\in S^2(V)\,\, | \, \ad w(\mu) = 0  \}.
$$
{\bf b.} The isomorphism $*: S^p(V) \to S^{m-p}(V)$ is equivariant with respect to the natural action of $\mathfrak{so}(V)$ on $S^*(V)$.
In particular,
$$
\mathrm{IDer} (\mu) = \mathrm{IDer} (*\mu).
$$

}

\noindent{\it Proof.} {\bf a.}  First of all using the Jacobi identity for $S^*V$ we obtain:
\begin{align*}
\ad w(\{ v_1,\dots, v_p\}) = [w, [v_1,\ldots,[v_n,\mu]\ldots]] =
\sum\limits^n_i [v_1,\ldots,[[w,v_i]\ldots,[v_n,\mu]]\ldots] +\\
[v_1,\ldots,[v_n,[w,\mu]]\ldots] = \sum\limits_j \{v_1,\ldots, [w,v_j],\ldots, v_n\} + [v_1,\ldots,[v_n,[w,\mu]]\ldots].
\end{align*}
We see that $\ad w$ is a derivation of $(V,\mu)$ if and only if $[w,\mu] = 0$.

{\bf b.} Let $L=e_1\dots e_m$ be as above and  $w\in S^2(V)$.
We have,
\begin{align*}
*([w,e_{i_1}  \dots e_{i_p}]) = *\Bigl(\sum\limits_{j=1}^p
 e_{i_1}  \dots [w,e_{i_j}]& \dots  e_{i_p}\Bigr)=\\
  & \sum\limits^p_{j=1} [e_{i_1},\ldots,[[w, e_{i_j}]\ldots,[e_{i_p},L]]\ldots].
 \end{align*}
 On the other side,
 $$
 \begin{array}{c}
[w, *(e_{i_1}  \dots e_{i_p})] = [w, [e_{i_1},  \dots[e_{i_p} , L]]]=
 \sum\limits^p_{j=1}  [e_{i_1},\ldots,[[w, e_{i_j}]\ldots,[e_{i_p},L]]\ldots].
\end{array}
$$
We use here the fact that $[w,L]=0$. Therefore, the $*$-operator is $\mathfrak{so}(V)$-equivariant.

Furthermore, assume that $w\in \mathrm{IDer} (\mu)$ or equivalently that $[w,\mu]=0$. Therefore,
$$
[w,*\mu]=*([w,\mu]) =*(0) =0.
$$
Hence, $w\in \mathrm{IDer}(*\mu)$. Conversely, if $w\in
\mathrm{IDer}(*\mu)$ then
$$
*([w,\mu]) = [w,*\mu]= 0.
$$
This finishes the proof. $\Box$

\subsection{Hodge decomposition for real metric strongly homotopy algebras}

\subsubsection{\bf Hodge decomposition for a vector space.}

In this Subsection we follow Kostant's approach \cite[Page 332 - 333]{Kos}.
 Let $W$ be a finite dimensional vector space with two linear operators $\d$ and $\delta$.

\medskip

\de\label{de Kostant disjoint operators}{\bf [Kostant]} Linear maps $\d$ and $\delta$ are called {\it disjoint} if the following holds:
\begin{enumerate}
\item $\d\circ \delta (x) = 0$ implies $\delta (x)=0$;

\item $\delta\circ \d (x) = 0$ implies $\d (x)=0$.
\end{enumerate}

\medskip

Denote $\mathcal{L} = \delta\circ \d + \d\circ \delta$.

\medskip

\prop{\bf [Kostant]} \label{Kostant Hodge decomposition} {\it Assume that $\d$ and $\delta$ are disjoint and $\d^2 = \delta^2 =0$. Then we have
$
 \Ker(\mathcal{L}) =  \Ker(\d)\cap \Ker(\delta)
$
and the direct sum (an analog of a Hodge Decomposition):
$$
W = \Im (\d) \oplus \Im (\delta) \oplus \Ker(\mathcal{L}).
$$
 In this case the restriction $\pi|_{\Ker(\mathcal{L})}$ of the canonical mapping }
$$
\pi : \Ker (\d) \to \Ker (\d)/\Im(\d) =: H(W,\d)
$$
{\it is a bijection. In other words $\Ker(\mathcal{L}) \simeq H(W,\d).\Box$}

\medskip

We will use this Proposition to obtain a Hodge decomposition for metric $L_{\infty}$-algebras.

\subsubsection{\bf Hodge decomposition for real metric  $L_{\infty}$-algebras.}

Let $V$ be a pure odd real $m$-dimensional vector space with a non-degenerate
skew-symmetric positive defined form $(\,,)$.
We can define a bilinear product $\langle\, ,\rangle$ in $S^* V$ by the following formula:
$$
\langle v_1 , v_2\rangle L = \left\{
\begin{array}{cl}
(-1)^{\frac{p(p-1)}{2}} v_1 \cdot (*v_2),\,\,& \text {if} \,\, v_1,v_2 \in S^p V;\\
0, \,\, & \rule{0pt}{5mm}\text {if $v_1\in S^p V$, $v_2 \in S^q V$ and $p\ne q$.}
\end{array}
\right.
$$
This bilinear product has the following properties:

\medskip

\prop\label{prop_inner product is sym and positive def} {\it Let us take $I = (i_1,\ldots, i_p)$ and $J = (j_1,\ldots, j_p)$ such that $i_1<\cdots< i_p$ and $j_1<\cdots< j_p$. Denote $e_{I}:= e_{i_1}\cdots e_{i_p}$ and $e_{j}:= e_{j_1}\cdots e_{j_p}$ We have
$$
\langle  e_{I}, e_{J} \rangle = \left\{
\begin{array}{c}
0, \,\, \text{if} \,\, I\ne J,\\
\rule{0pt}{4mm}1, \,\, \text{if} \,\, I= J.
\end{array}
\right.
$$
In particular, the pairing $\langle\, ,\rangle$ is symmetric and positive definite.}
\medskip

\noindent{\it Proof.} A straightforward computation.$\Box$

\medskip

Let $\mu \in S^*V$ be
 any element. Denote by $\d:S^*V\to S^*V$ the linear operator  $v
 \mapsto [\mu,v]$.  Let $\mu= \sum\limits_k\mu_k$, where $\mu_k\in S^{k+2}V$,  and we put $\d_k:= [\mu_k, -]$. 
 Using Hodge
 $*$-operator we can define the following operator
 $$
 \delta = \sum\limits_k (-1)^{\frac{k(1-k)}{2}}\delta_k,
 $$
 where $\delta_k:=*\d_k *$.

\medskip

\prop\label{prop d and delta are adoint}{\it Assume that $\mu \in  S^*(V)$,  $\d$ and $\delta$ are as above. Then we have
$$
\langle \d(v) , w\rangle = - (-1)^{\frac{m(m-1)}{2}} \langle v , \delta(w)\rangle
$$
 for $v,w \in S^{*} V$.  The operators $\d$ and $\delta$ are disjoint.}

\medskip

\noindent{\it Proof.}  Let us take $\mu_k\in S^{k+2}
V$, $v \in S^{p-k} V$ and $w \in S^p V$.  (We assume that $S^{r} V
=\{0\}$ for $r<0$ and $r>m$,  where $m = \dim V$.) Then,  $v\cdot
*w\in S^{m-k}V$ and we have:
$$
[\mu_k,v\cdot *w] \subset  [S^{k+2} V,S^{m-k}(V)] =0.
$$
Furthermore,
\begin{align*}
0= [\mu_k, v\cdot *w] = [\mu_k, v] \cdot *w  + (-1)^{\bar{\mu}_k \bar{v}} v\cdot [\mu_k, *w] =
[\mu_k&, v] \cdot *w 
+ \\
 (-1)^{\bar{\mu}_k\bar{v} + \frac{m(m-1)}{2}} v\cdot ** [\mu_k, *w] = 
\d_k (v) \cdot *w +& (-1)^{\bar{\mu}_k\bar{v} + \frac{m(m-1)}{2}} v\cdot * \delta_k (w),
\end{align*}
where $\d_k(v) = [\mu_k, v]$ and $\delta_k(w) =  *[\mu_k, *w]$. Further,
\begin{align*}
\d_k (v) \cdot *w = (-1)^{\frac{p(p-1)}{2}} \langle \d_k (v), w \rangle L;\quad
v\cdot * \delta_k (w) = (-1)^{\frac{(p-k)(p-k-1)}{2}} \langle v, \delta_k (w)\rangle L.
\end{align*}
Therefore,
$$
(-1)^{\frac{p(p-1)}{2}}\langle\d_k(v), w \rangle  = - (-1)^{\frac{(p-k)(p-k-1)}{2}} (-1)^{\bar{\mu}_k\bar{v} + \frac{m(m-1)}{2}} \langle v,\delta_k (w)\rangle
$$
or 
$$
\langle\d_k(v), w \rangle  = - (-1)^{\frac{k(1-k)}{2}} (-1)^{ \frac{m(m-1)}{2}} \langle v,\delta_k (w)\rangle
$$
for all $v \in S^{p-k} V$ and $w \in S^p V$. 
Note that this equation holds trivially for $v\in S^s V$ and $w\in S^q V$, where $q-s\ne k$. 
Therefore, we have
\begin{equation}\label{eq (d_kv,w)=+-(v,delta_k w)}
\langle\d_k(v), w \rangle  = - (-1)^{\frac{k(1-k)}{2}} (-1)^{ \frac{m(m-1)}{2}} \langle v,\delta_k (w)\rangle
\end{equation}
for all $v,w \in S^{*} V$ and $\mu_k\in S^{k+2} V$.

 Let us take any $\mu \in S^*(V)$. Then $\mu = \sum\limits_k \mu_k$, where $\mu_k\in  S^k(V)$. Therefore, $\d$ and $\delta$ also possess corresponding decomposition:
   $\d= \sum\limits_k \d_k$ and $\delta = \sum\limits_k (-1)^{\frac{k(1-k)}{2}} \delta_k$, where $\d_k= [\mu_k,-]$ and $\delta_k =  *\d_k*$.
Using (\ref{eq (d_kv,w)=+-(v,delta_k w)}), we get for any $v,w \in S^{*} V$:
\begin{align*}
\langle\d(v), w \rangle = \sum\limits_k \langle\d_k(v), w \rangle =  - \sum\limits_k (-1)^{\frac{k(1-k)}{2}} (-1)^{ \frac{m(m-1)}{2}} \langle v,\delta_k (w)\rangle =\\
- (-1)^{ \frac{m(m-1)}{2}}\sum\limits_k (-1)^{\frac{k(1-k)}{2}}\langle v,\delta_k (w)\rangle =
 - (-1)^{ \frac{m(m-1)}{2}} \langle v,\delta (w)\rangle.
\end{align*}
 The first statement is proven.

 Let us show that $\d\circ \delta (v) = 0$ implies
$\delta (v) = 0$, i.e. the  operators $\d$ and $\delta$ are
disjoint. (This argument we borrow from \cite{Kos}.) Indeed,
$$
0= \langle\d\circ \delta (v), v\rangle = - (-1)^{ \frac{m(m-1)}{2}} \langle\delta (v), \delta (v)\rangle.
$$
The pairing $\langle\,,\rangle$ is positive definite, hence $\delta (v) =0$. Analogously we can show that $\delta\circ \d (v) = 0$ implies $\d (v) = 0$.$\Box$

\medskip

Assume that $(V,\mu)$ is an $L_{\infty}$-algebra. By Proposition \ref{prop quadratic L_infinity algebra} this means that $\mu$ is an odd element and $[\mu,\mu]\in \mathbb K$.  Denote by $H(V,\mu)$ the cohomology group of $(V,\mu)$. By definition $H(V,\mu) := \Ker
(\d)/ \Im (\d)$, where $\d=[\mu,-]:S^*V \to S^*V$. (Clearly, $\d^2= \delta^2=0$.) The main result of this section is the following
theorem.

\medskip

\t\label{Teor Hodge decomposition for real Lie algebras} {\bf [Hodge
decomposition  for real metric $L_{\infty}$-algebras]}
{\it Let $\mu \in S^*(V)$ be a real metric $L_{\infty}$-algebra
structure on $V$ and  $\d$ and $\delta$ be as above. Then we have a
direct sum decomposition:
$$
V = \Im (\d) \oplus \Im (\delta) \oplus \Ker(\mathcal{L}),
$$
where $\mathcal{L} = \delta \circ \d + \d\circ \delta$, and $
 \Ker(\mathcal{L}) \simeq H(V,\mu).
$}

\medskip

\noindent{\it Proof.} The statement follows from Propositions \ref{Kostant Hodge decomposition} and \ref{prop d and delta are adoint} $\Box$.

\section{Filippov and Lie $m$-dimensional invariant\\ $(m-3)$-algebras }

\subsection{Anti-commutative $m$-dimensional invariant $(m-3)$-ary algebras and coadjoint orbits}

In this section  we will
classify all $(m-3)$-ary anti-commutative algebra structures on $V$, where $\dim V=m$, up to orthogonal
isomorphism in terms of coadjoint orbits. Let
again $V$ be a pure odd vector space with an even non-degenerate
skew-symmetric form $(\,,)$. That is $V=V_{\bar 1}$ and $(a,b) = (b,a)$ for all $a,b\in
V$. As usual we denote by $\mathrm{O}(V)$ the Lie group of all
invertible linear operators on $V$ that preserve the form $(\,,)$
and by $\mathrm{SO}(V)$ the subgroup of $\mathrm{O}(V)$ that
contains all operators with the determinant $+1$. We have
$\mathfrak{so}(V)= \Lie \mathrm{O}(V)= \Lie \mathrm{SO}(V)$.

\medskip

\de Two $n$-ary algebra structures $\mu,\mu'\in S^* V$ on $V$ are called {\it isomorphic} if there exists $\varphi \in \mathrm{SO}(V)$ such that
$$
\varphi(\{v_1,\ldots, v_n\}_{\mu}) = \{\varphi(v_1),\ldots,\varphi( v_n)\}_{\mu'}
$$
for all $v_i\in V$. Here we denote by $\{\ldots\}_{\nu}$ the multiplication on $V$ corresponding to the algebra structure  $\nu$.

\medskip

Sometimes we will consider isomorphism of $n$-ary algebra structures up to $\varphi \in \mathrm{O}(V)$.
We need the following two lemmas:

\medskip

\lem \label{lem O preserves brackets} {\it Let us take $\varphi \in\mathrm{O}(V)$ and $w,v\in S^*V$.  Then, $\varphi([w,v]) = [\varphi(w),\varphi(v)]$.}

\medskip

\noindent{\it Proof.} It follows from the following two facts:
\begin{itemize}
\item $(\varphi(w),\varphi(v))= (w,v)$,\, if $w,v\in V$;

\item $\varphi(w\cdot v) = \varphi(w)\cdot \varphi( v)$ for all $w,v\in S^*V$. $\Box$

\end{itemize}

\medskip

\lem \label{lem mu i mu' are isomorphic} {\it  Two $n$-ary  algebra structures
$\mu$ and $\mu'$, where $\mu, \mu'\in S^{n+1}V$, are isomorphic
if and only if there exists $\varphi \in \mathrm{SO}(V)$ such that
$\varphi(\mu) = \mu'$. In other words, two  $n$-ary algebra structures are
isomorphic if and only if they are in the same orbit of the action $\mathrm{SO}(V)$ on
$S^{n+1}V$.}

\medskip

\noindent{\it Proof.} From Lemma \ref{lem O preserves brackets} it follows that
$$
\varphi(\{v_1,\ldots, v_n\}_{\mu}) = \{\varphi(v_1),\ldots,\varphi( v_n)\}_{\varphi(\mu)}.
$$
Therefore, $n$-ary algebra structures $\mu$ and $\varphi(\mu)$ are isomorphic. 

Conversely, if $\mu$ and $\mu'$ are isomorphic and $\varphi
\in \mathrm{SO}(V)$ is an isomorphism then from the definition it
follows that:
$$
\{\varphi(v_1),\ldots,\varphi( v_n)\}_{\varphi(\mu)} = \{\varphi(v_1),\ldots,\varphi( v_n)\}_{\mu'}
$$
for all $v_i\in V$. Therefore, $\varphi(\mu) =\mu' $.$\Box$

\medskip

\medskip
\t \label{theor classification of m-3 ary algebras coad orbits} {\it Assume that $\dim V = m$.
Classes of isomorphic real or complex invariant $(m-3)$-ary  algebra structures on $V$ are in one-to-one correspondence with
coadjoint orbits of the Lie group $\mathrm{SO}(V)$. }

\medskip

\noindent{\it Proof.} Note that
in the case of the Lie group $\mathrm{SO}(V)$ the adjoint and coadjoint actions are equivalent. By Proposition \ref{prop derivations},b, the isomorphism $*: S^2(V)\to S^{m-2}(V)$ is $\mathfrak{so}(V)$-equivariant. Clearly, it is also $\mathrm{SO}(V)$-equivariant and the action of $\mathrm{SO}(V)$ on $S^2(V)$ is equivalent to the adjoint action of $\mathrm{SO}(V)$. (Note that $S^2(V)\simeq \mathfrak{so}(V)$ by Proposition \ref{prop lin operators preserving the form}.) The result follows from 
Lemma \ref{lem mu i mu' are isomorphic}.$\Box$

\medskip

\ex\label{ex classification m dim (m-2)-ary} Assume that $m \geq 4$. Any $m$-dimensional invariant real or complex $(m-2)$-ary algebra has the form $(V,\mu)$, where $\mu = *(v)$ and $v\in V$. This follows from the existence of the isomorphism $*:V\to \bigwedge^{m-1}V$. 

\medskip

Assume that $\mathbb K=\mathbb R$. It is well-known that any real skew-symmetric matrix $A$ can be written in the form $A= QA'Q^{-1},$
where
$$
A'= \mathrm{diag}(J_{a_1}, \ldots, J_{a_k}, 0),\quad m=2k+1,
$$
and
$$
A'= \mathrm{diag}(J_{a_1}, \ldots, J_{a_k}),\quad m=2k.
$$
Here
$$
J_{a_j} = \left(
\begin{array}{cc}
0& a_j\\
-a_j& 0
\end{array}
 \right), \quad a_j\in \mathbb{R},
$$
and $Q\in \mathrm{SO}(V)$. Moreover we can assume that $0 \leq a_k \leq \cdots \leq a_1$ if $m$ is odd and $ |a_k| \leq \cdots \leq a_1$ if $m$ is even. Therefore, coadjoint orbits are parametrized by the numbers $(a_j)$.

Let $(\xi_i)$ be an orthogonal basis
of $V$ such that the matrix $A\in S^2(V) \simeq  \mathfrak{so}(V)$ has the form
$$
A=\mathrm{diag}(J_{a_1}, \ldots, J_{a_k}, 0) \,\,\,\text{or}\,\,\, A=\mathrm{diag}(J_{a_1}, \ldots, J_{a_k}).
$$
Then the corresponding element in $S^2V$ is
$$
v_A= a_1 \xi_1\xi_2 + \ldots, a_k \xi_{2k-1}\xi_{2k},
$$
where $a_j$ are as above. We obtained the following theorem:

\medskip

\t \label{theor classification of real m-3 ary
algebras}{\bf[Classification of real invariant $(m-3)$-ary $m$-dimensional algebras] }{\it Real invariant $(m-3)$-ary $m$-dimensional algebras are parame\-trized by vectors
$$
v= a_1 \xi_1\xi_2 + \ldots, a_k \xi_{2k-1}\xi_{2k},
$$
where $0 \leq a_k \leq \cdots \leq a_1$ if $m=2k+1$ and $ |a_k| \leq \cdots \leq a_1$ if $m=2k$.
Explicitly such algebras are given by the derived potentials $\mu_v$, where $\mu_v= *(v).$

}

\medskip

\subsection{Classification of real $m$-dimensional simple $(m-3)$-ary algebras with a positive definite invariant form}

In this section we give a classification of real $m$-dimensional simple $(m-3)$-ary algebras with an invariant form up to orthogonal isomorphism. Let
 $V$ be a pure odd vector space over $\mathbb R$ with an even non-degenerate 
skew-symmetric form $(\,,)$. That is $V=V_{\bar 1}$ and $(a,b) = (b,a)$ for all $a,b\in V$. We assume in addition that $(\,,)$ is positive definite. 

\medskip

\de A vector subspace $W\subset V$ is called an {\it ideal} of a symmetric $n$-ary algebra $(V,\mu)$ if $\mu(V, \ldots, V, W) \subset W$.

\medskip

In other words, the vector space $W$ is an ideal if and only if it
is invariant with respect  to the set of endomorphisms $\mu(v_1,
\ldots, v_{n-1},-): V\to V$, where $v_i\in V$. Clearly, the vector
space $W$ is an ideal if and only if it is invariant with respect to
the Lie algebra $\mathfrak{g}$ generated by all $ \mu(v_1,\ldots, v_{n-1},-)$.

\medskip

\de A symmetric $n$-ary algebra $(V,\mu)$ is called {\it simple} if it is not $1$-dimensional and if it does not have any proper ideals.

\medskip

\ex\label{ex Ling classification} The classification of simple complex and real Filippov  $n$-algebras  was done in \cite{Ling}: there is one series of complex
Filippov $n$-algebras $A_k$, where $k$ is a natural number and
several real forms for each $A_k$. The complex $n$-ary algebra $A_k$ has an invariant form and in our terminology it is given by the derived potential $L$ (a top form) and formula (\ref{eq der product}).

\medskip

\ex\label{ex classification m=5} Let $m=5$. By Theorem \ref{theor classification of real m-3 ary algebras} we see that we have three types of $2$-ary algebras up to isomorphism ($\mu_i$ is a derived potential):
\begin{itemize}
\item $\mu_1 = *(0)=0$;
\item $\mu_2 = *(a_1\xi_1\xi_2) = b_1 \xi_3\xi_4\xi_5$, where $b_1= \pm a_1\ne 0$;
\item $\mu_3 = *(a_1\xi_1\xi_2+ a_2\xi_3\xi_4)= b_1 \xi_3\xi_4\xi_5 + b_2 \xi_1\xi_2\xi_5$, where $b_1= \pm a_1\ne 0$ and $b_2= \pm a_2\ne 0$.
\end{itemize}
Obviously, the zero algebra $(V,0)$ is not simple. The derived potential $\mu_2 = b_1 \xi_3\xi_4\xi_5$ corresponds to the algebra with a non-trivial center:
$$
(\mu_2)_{\xi_1} = (\mu_2)_{\xi_2} =0.
$$
Therefore, $[x,[\xi_i,\mu]] = 0$, $i=1,2$, for any $x\in V$. Hence $\langle\xi_1,\xi_2 \rangle$ is an ideal and the algebra $(V,\mu_2)$ is not simple. We will see that the algebra
$\mu_3$ is simple. It is not a  Lie algebra because $[\mu_3,\mu_3]=
-2b_1b_2 \xi_1\xi_2\xi_3\xi_4 \ne 0$.

\bigskip

\t \label{theor classification of real simple m-3 ary
algebras}{\bf[Classification of real simple $(m-3)$-ary algebras
with an invariant form] }{\it Assume that $m>4$. All real $m$-dimensional $(m-3)$-ary
algebras from Theorem \ref{theor classification of real m-3 ary algebras} are simple except of two cases:

\begin{itemize}
\item $v=0$;

\item $v=  a_1 \xi_1\xi_2$, where $a_1\ne 0$.
\end{itemize}

}

\medskip

\noindent{\it Proof.} As in the $2$-ary case we can show that if $W$ is an ideal in $V$, then $W^{\perp}$ is also an ideal in $V$. Let us take any real $m$-dimensional $(m-3)$-ary
algebra $(V,\mu)$ and let us assume that $W$ is an ideal. Then $V=W\oplus W^{\perp}$ and we have the decomposition:
$$
\mu = \sum_{t}\mu_t, \quad \text{where} \quad \mu_t\in \bigwedge^{t}W \wedge \bigwedge^{m-2-t}W^{\perp}.
$$
Since $W$ and $W^{\perp}$ are ideals, we have $\mu_s=0$ for $s\ne 0, m-2$.  Assume that $\mu\ne 0$, then one of these ideals has dimension greater than or equal to $m-2$. Hence, we can assume that $\dim W = 1$ or $2$. In case $\dim W =1$, we have $\mu_{m-2}=0$ and $(W^{\perp}, \mu_{0})$ is an $(m-1)$-dimensional $(m-3)$-ary algebra, where $\mu_0= *(w)$ and $w\in W^{\perp}$, see Example \ref{ex classification m dim (m-2)-ary}. The algebra $(W^{\perp}, \mu_{0})$ has a zero ideal $\langle w\rangle$, since $[w,\mu_0]=0$. Therefore, we can assume that $\dim W=2$. In case $\dim W=2$, we have $\mu_{m-2}=0$ and $(W^{\perp}, \mu_{0})$ is an $(m-2)$-dimensional $(m-3)$-ary algebra. Hence, $\mu_{0}\in \bigwedge^{top}W^{\perp}$. If $\mu_0\ne 0$, then the algebra $(W^{\perp}, \mu_{0})$ is simple, see Example \ref{ex Ling classification}. Therefore, any real invariant $(m-3)$-ary algebra with a proper ideal has the form $(V,0)$ or $(V,*(a_1 \xi_1\xi_2))$. All other algebras are simple.$\Box$

\bigskip

\subsection{ Classification of real $m$-dimensional invariant simple $(m-3)$-ary algebras 
 satisfying Jacobi identity \ref{eq_ Jacobi Filippov} and \ref{eq_ Jacobi SH}}

In this Section we classify real simple $n$-ary algebras with a positive definite invariant form satisfying Jacobi identity \ref{eq_ Jacobi Filippov} and \ref{eq_ Jacobi SH}. 

\medskip

\noindent {\bf Jacobi identity \ref{eq_ Jacobi Filippov}.}  In \cite{Ling} it was proven that there exist only one complex Filippov $n$-algebra for any $n>2$. This algebra is $(n+1)$-dimensional. In our notations it is given by $*(1)=L$. Another result in \cite{Ling} is the following: 

\smallskip
{\it A real simple Filippov $n$-algebra is isomorphic to the realification of a simple complex Filippov $n$-algebra or to a real form of a simple complex Filippov $n$-algebra.}

\smallskip
In particular real simple Filippov $n$-algebras are of dimension $n+1$ or $2n+2$. It follows that simple $n$-ary algebras in Theorem 
\ref{theor classification of real simple m-3 ary
algebras} are not of Filippov type. 
For $n=m-2$ any derived potential has the form $\mu = *(v)$, where
$v\in V\setminus \{0\}$. All such algebras have non-trivial centers
because $[v,\mu]=0$. Therefore, they are not simple.
Furthermore, such algebras are of Filippov type. Indeed,  since $L$ satisfy the Jacobi identity (\ref{eq_ Jacobi Filippov}) by Proposition \ref{prop quadratic Filippov algebra}, we have $[L_{a_1,\dots, a_{m-1}}, L] = 0$ for any $a_i\in V$. Furthermore, for $\mu = *(v)=[v,L]$ we get
$$
[\mu_{a_1,\dots, a_{m-2}},\mu] = [L_{a_1,\dots, a_{m-2},v}, L_v ] =  [v, [L_{a_1,\dots, a_{m-2},v}, L] ] = 0.
$$  
By Proposition \ref{prop quadratic Filippov algebra}, we see that $(V,\mu)$ is a Filippov algebra. By the same argument the derived potential $[v,[w,L]]$ also corresponds to a Filippov algebra. 

\medskip

\t \label{teor classification of real Filippov} {\it  Assume that $m>4$. Real $m$-dimensional Filippov $n$-algebras with a positive definite invariant form, where $n=m-1$, $m-2$ or $m-3$, are given up to $\mathrm{SO}(V)$-isometry by the following derived potentials:
\begin{itemize}
\item $\mu = *(0) =  0$;

\item $\mu = *(a \cdot 1)= a\xi_1\cdots \xi_{m}$, where $a\in \mathbb R \setminus \{0\}$;  

\item $\mu = *(a\xi_1) = a\xi_2\cdots \xi_{m}$, where $a\in \mathbb R^{>0}$; 

\item $\mu = *(a \xi_1\xi_2) = -a\xi_3\cdots \xi_{m}$, where $a\in \mathbb  R^{>0}$.

\end{itemize}
\medskip
}

\noindent {\bf Jacobi identity (\ref{eq_ Jacobi SH}).} 
As above assume that $m>4$ and $(\,,)$ is a symmetric positive definite form. 

\medskip

\t \label{teor classification of real sh} {\it All 
algebras in Theorem \ref{theor classification of real m-3 ary algebras} satisfy Jacobi identity (\ref{eq_ Jacobi SH}) with the exception of the following cases:
\begin{itemize}
\item $m=5$, the algebras with the derived potential $\mu = *(a_1\xi_1\xi_2 + a_2 \xi_3\xi_4)$, where $a_1,a_2\ne 0$; 
\item $m=6$, the algebras with the derived potentials $\mu = *(a_1\xi_1\xi_2 + a_2 \xi_3\xi_4)$ and $\mu = *(a_1\xi_1\xi_2 + a_2 \xi_3\xi_4 + a_3 \xi_5\xi_6)$, where $a_i\ne 0$; 
\end{itemize}
}

\medskip

\noindent {\it Proof.} Assume that  $m$ is odd. By Corollary of Proposition \ref{prop quadratic L_infinity algebra} in this case  Jacobi identity \ref{eq_ Jacobi SH} is equivalent to $[\mu,\mu]=0$. Assume that $m>5$, then $[\mu,\mu]\in S^{2m-6}V = \{0\}$. In the case $m=5$ the result follows from  Example \ref{ex classification m=5}.  

Assume that $m$ is even. First of all consider the case $m=6$. Let us take 
$$
\mu = b_1 \xi_3\xi_4 \xi_5 \xi_6 + b_2 \xi_1\xi_2 \xi_5 \xi_6, \,\,\, b_1,b_2\ne 0.
$$
Denote by $LHS$ the left hand side of (\ref{eq_ Jacobi SH}). Let us calculate $LHS$ for $(\xi_1,\ldots,\xi_5)$. 
$$
LHS= \{ \{\xi_1,\xi_2,\xi_5 \}, \xi_3,\xi_4 \} + \{ \{\xi_3,\xi_4,\xi_5 \}, \xi_1,\xi_2 \} = -2b_1b_2\xi_5\ne 0.
$$
The main idea here is to use the fact that $\{x,y,z\}=0$ if $x\in \{\xi_1,\xi_2\}$ and  $y\in \{\xi_3,\xi_4\}$. The proof for 
$$
\mu = b_1 \xi_3\xi_4 \xi_5 \xi_6 + b_2 \xi_1\xi_2 \xi_5 \xi_6 + b_3 \xi_1\xi_2 \xi_3 \xi_4, \,\,\, b_i\ne 0
$$
is similar.

Consider the case $m$ is even and $m>6$. Let us compute $LHS$ for $(v_i)$, where $i=1,\ldots, 2m-7$. Without loss of generality we can assume that between elements $v_i$ are at least two equal. Let $v_s=v_t=v$. Clearly, $\{v_{i_1},\ldots,v,\ldots,v,\ldots, v_{i_n}\} = 0$. 
Therefore,
$$
LHS = \sum\limits_{k,l} J_1^{(k,l)} + \sum\limits_{k,l} J_2^{(k,l)},
$$
where $J_1^{(k,l)}$ and $J_2^{(k,l)}$ is the sum of all summands of the form  
$$
\begin{array}{c}
\{\{v_{i_1},\ldots,\underset{k}{v_s},\ldots, v_{i_{m-3}}\}, v_{j_1},\ldots,\underset{l}{v_t},\ldots, v_{j_{m-4}} \},\\

\{\{v_{i_1},\ldots,\underset{k}{v_t},\ldots, v_{i_{m-3}}\}, v_{j_1},\ldots,\underset{l}{v_s},\ldots, v_{j_{m-4}} \}
\end{array}
$$
respectively.
Further,
$$
\begin{array}{rl}
J_1^{(k,l)} =  &\pm \sum (-1)^{(I,J)} \{\{v_{i_1},\ldots,\underset{k}{\hat{v}_s},\ldots, v_{i_{m-3}}, v_s\}, v_{j_1},\ldots,\underset{l}{\hat{v}_t},\ldots, v_{j_{m-4}}, v_t \} = \\

&\pm \sum (-1)^{(I',J')} \{\{v_{i_1},\ldots,\underset{k}{\hat{v}_s},\ldots, v_{i_{m-3}}\}_{v}, v_{j_1},\ldots,\underset{l}{\hat{v}_t},\ldots, v_{j_{m-4}} \}_v,
\end{array}
$$
where $\{\ldots \}_{v}$ is the multiplication corresponding to the derived potential $\mu_v = [v,\mu]$ and $(-1)^{(I',J')}$ is the sign of the permutation 
$$
 (v_1, \ldots,\hat{v}_s, \ldots, \hat{v}_t, \ldots, v_{2m-7})\longmapsto (v_{i_1},\ldots,\underset{k}{\hat{v}_s},\ldots, v_{i_{m-3}}, v_{j_1},\ldots,\underset{l}{\hat{v}_t},\ldots, v_{j_{m-4}} ).
 $$
Since $\mu_v\in S^{m-3} W$, where $W = \langle v\rangle^{\bot}$, we see that $[\mu_v,\mu_v] = 0$. Therefore (\ref{eq_ Jacobi SH}) holds for $\{\ldots \}_{v}$ and 
$J_1^{(k,l)} = 0$. Similarly, $J_2^{(k,l)} = 0$. The proof is complete.$\Box$

\medskip

\noindent {\bf Corollary.}  {\it All simple algebras from Theorem \ref{theor classification of real simple m-3 ary algebras} satisfy Jacobi identity (\ref{eq_ Jacobi SH}) for $m>6$. }

\section{Quasi-Frobenius skew-symmetric $n$-ary\\
	 algebras}

Let $V$ be a pure odd vector space and $\mu \in S^{n}(V^*)\otimes V$
be an $n$-ary  symmetric algebra structure on $V$.

\medskip

 \de \label{de quasi-Fronebius algebras} An $n$-ary algebra $(V,\mu)$ is called {\it quasi-Frobenius}
 if it is equipped with a symmetric bilinear form $\varphi$ such that
\begin{equation}\label{eq quasi-Fronebius algebras}
\sum_{cycl} \varphi(a_1,\mu(a_2,\ldots, a_{n+1}))  = 0.
\end{equation}

\smallskip

If we forget about superlanguage this means that the algebra
$(V,\mu)$ is skew-symmetric and $\varphi$ is a
 skew-symmetric bilinear form on $V$.

\medskip

\ex Assume that $n=2$ and $(V,\mu)$ is a Lie algebra. Then our
definition coincides with the definition of a  quasi-Frobenius Lie
algebra. Recall that a {\it quasi-Frobenius Lie algebra} is a Lie
algebra $\mathfrak g$ equipped with a non-degenerate skew-symmetric
bilinear form $\beta$ such that
$$
\beta([x,y],z) + \beta([z,x],y) +\beta([y,z],x) =0.
$$

\medskip

We may assign an $n$-ary algebra $(V\oplus V^*, \mu^T)$ to
$(V,\mu)$, called the  $T^*_0$-extension of $(V,\mu)$. (The notion
of  $T^*_{\theta}$-extension for algebras was introduced and studied
in \cite{Bordemann}. We will need this notion only for $\theta =0$.)
The construction of  $(V\oplus V^*, \mu^T)$ is the following: the
$n$-ary algebra structure $\mu^T$ is just the image of $\mu$ by the
natural inclusion $S^{n}(V^*)\otimes V \hookrightarrow S^{*}(V^*
\oplus V)$. Furthermore, the pure odd vector space $ V\oplus V^*$
has a skew-symmetric (in supersense) pairing given by
$$
(a,\alpha) = (\alpha, a)= \alpha(a),
$$ where $\alpha\in V^*$ and $a\in V$. This defines a Poisson bracket on $S^*(V\oplus V^*)$. So
 $(V\oplus V^*, \mu^T)$ as a quadratic symmetric $n$-ary algebra, where
 the multiplication is given by the derived bracket with the derived potential $\mu^T\in S^{*}(V^* \oplus V)$.
 More precisely,  the new multiplication  $\mu^T$ in $ V\oplus V^*$ is given by:
$$
\mu^T|_{S^n(V)} = \mu,\quad \mu^T|_{S^{n-k}(V)\cdot S^{k}(V^*)} = 0 \,\,\,\text{if}\,\,\, k>1,\quad \mu^T(S^{n-1}(V)\cdot S^{1}(V^*)) \subset V^*
$$
and
$$
\mu^T(a_1,\ldots,a_{n-1}, b^*) (c) := - b^*(\mu(a_1,\ldots,a_{n-1}, c)).
$$

The main observation here is:

\medskip

\prop\label{prop quasi frobenius algebras} {\it Let
 $V$ be a pure odd vector space and $n$ be even.
 Then an $n$-ary algebra $(V,\mu)$ has a  quasi-Frobenius structure with
 respect to a symmetric  form $\varphi$ if and only if the maximal isotropic
 subspace $B_{\varphi} = \{ a+ \varphi(a,-) \} \subset V\oplus V^*$ is a subalgebra in $( V\oplus V^*,\mu^T)$.
In other words, there is a one-to-one correspondence between quasi-Frobenius structures
on $(V,\mu)$ and maximal isotropic subalgebras in  $( V\oplus V^*,\mu^T)$ that are transversal to $V^*$.}

\smallskip

\noindent {\it Proof.} First of all it is well-known that maximal isotropic subspaces in
 $V\oplus V^*$ that are transversal to $V^*$ are in one-to-one correspondence with
 $\varphi\in S^2V$. Let us show that $\varphi$ satisfies (\ref{eq quasi-Fronebius algebras})
 if and only if $B_{\varphi}$ is a subalgebra. Denote $a^*:= \varphi(a,-)\in V^*$. Then we have:
\begin{align*}
(\mu^T(a_1+ a_1^*, \ldots,a_n&+ a_n^*), c+ c^*) =\\ &c^*(\mu(a_1,\ldots, a_n)) + \sum\limits_k(\mu^T(a_1,\ldots,a^*_k,\ldots,a_n), c)= \\
&\varphi(c,\mu(a_1,\ldots, a_n)) - \sum\limits_k a^*_k (\mu(a_1,\ldots, a_{k-1},c, a_{k+1},\ldots,a_n)) =\\
&\varphi(c,\mu(a_1,\ldots, a_n)) - \sum\limits_k \varphi(a_k, \mu(a_1,\ldots, a_{k-1},c, a_{k+1}\ldots,a_n)).
\end{align*}
Furthermore,
\begin{align*}
\varphi(a_k, \mu(a_{k+1},\ldots, a_n,c,a_1\ldots, a_{k-1})) =&\\ (-1)^{(k-1)(n-k-1) } & \varphi(a_k, \mu(a_{1},\ldots, a_{k-1},a_{k+1}\ldots,a_{n},c)) =\\
(-1)^{k(n-k-1) +1 } &\varphi(a_k, \mu(a_{1},\ldots, a_{k-1},c,a_{k+1}\ldots,a_{n})).
\end{align*}
If $n$ is even, $(-1)^{k(n-k-1) +1 }  = -1$. Therefore,  we have:
$$
(\mu^T(a_1+ a_1^*, \ldots,a_n+ a_n^*), a_{n+1}+ a_{n+1}^*) = \sum_{cycl} \varphi(a_1,\mu(a_2,\ldots, a_{n+1})).
$$
This expression is equal to $0$ if and only if the algebra $(V,\mu)$ is quasi-Frobenius with respect to $\varphi$.
On other side, $(\mu^T(a_1+ a_1^*, \ldots,a_n+ a_n^*), a_{n+1}+ a_{n+1}^*)$ is equal to $0$ if and only if $B_{\varphi}$ is a subalgebra in $( V\oplus V^*,\mu^T)$.
The proof is complete.$\Box$

\medskip

\noindent {\bf Remark.} The result of Proposition \ref{prop quasi frobenius algebras} is well-known for Lie algebras.

\noindent{\it Elizaveta Vishnyakova}

\noindent {Max Planck Institute for Mathematics Bonn and}

\noindent University of Luxembourg

 \noindent{\emph{E-mail address:}
\verb"VishnyakovaE@googlemail.com"}

\end{document}